\magnification=\magstep1
\documentstyle{amsppt}
\parindent=0pt
\TagsOnLeft
\advance\hoffset by 0.1 truein

\vsize=7.2in

\font\headlinerm=cmr10 scaled \magstep2

\font\smallrm=cmr8 
 
\font\largebf=cmbx12

\nopagenumbers
\document

\

\headlinerm Global Stability for Holomorphic Foliations in Kaehler Manifolds

\vskip0.4cm

\smallrm Jorge Vit\'{o}rio Pereira\footnote {Supported
by IMPA-CNPq}

\vskip0.4cm

\bf Abstract. \smallrm \ \ \ \ 
We prove the following theorem for Holomorphic Foliations in
compact complex kaehler manifolds: if there is a compact leaf with finite holonomy, then
every leaf is compact with finite holonomy. As corollary we reobtain
stability theorems for compact foliations in Kaehler manifolds of Edwards-Millett-Sullivan and
Hollman. 

\vskip0.4cm

\largebf 1. Introduction

\vskip0.4cm
\rm

The question of global stability is recurrent in the theory of foliations.
The work of Ehresmann and Reeb establishes the so called global stability
theorem, which says that if $\Cal F$ is a transversely orientable
codimension one foliation in a compact connected manifold $M$ that has a
compact leaf $L$ with finite fundamental group, then every leaf of $L$ is
compact with finite holonomy group[Ca-LN]. Counterexamples for codimension
greater than one are known since the birth of the theorem. Here we want to
abolish the hypothesis on the codimension for a special kind of foliation,
namely holomorphic foliations in complex Kaehler manifolds. In other words
we are going to prove the following :

\proclaim{Theorem 1}
{Let }$\Cal {F}${\ be a holomorphic foliation of
codimension q in a compact complex Kaehler manifold. If }$\Cal {F}$%
{\ has a compact leaf with finite holonomy group then every leaf of }$%
\Cal {F}${\ is compact with finite holonomy group.}
\endproclaim

\rm
Another kind of stability problem was posed by Reeb and Haefliger. The
question was the stability of compact foliations, that is, if a foliation
has all leaves compact is the leaf space Hausdorff? Positive answers to this
problem arose in the work of Epstein[Ep], Edwards-Millet-Sullivan[EMS], Holmann[Ho], etc.
There are plenty situations where the leaf space is not Hausdorff.
Sullivan found a example in the $C^\infty $ case[Su], Thurston in the
analytic case[Su] and M\"{u}ller in the holomorphic case[Ho]. The examples
of Sullivan and Thurston live in compact manifolds, and M\"{u}ller's in a non-compact non-Kaehler manifold. As corollary of the
theorem we reobtain Holmann's result and a special case of [EMS]'s outstanding Theorem.  

\proclaim{Corollary}[EMS,Ho]
{Suppose M is a complex Kaehler manifold. If }$\Cal {F}$%
{\ is a compact foliation, i.e., every leaf is compact, then every
leaf has finite holonomy group. Consequently, there is an upper bound on the
volume of the leaves, and the leaf space is Hausdorff.}
\endproclaim

The author would like to thanks B. Sc\'{a}rdua for valuable conversations.

\vskip0.8cm
\pagebreak

\largebf 2. The Leaf Volume Function

\vskip0.4cm

\rm

Let $\Cal {F}$ be a holomorphic foliation of a complex Kaehler manifold $%
(M,\omega )$. As in [Br] we define 
$$
\Omega =\{p\in M|\text{ the leaf }L_p\text{ through p is compact with finite
holonomy\}} 
$$

By the local stability theorem of Reeb[Ca-LN] $\Omega $ is an open set of $M 
$. Set, for every $p\in \Omega ,$ $n(p)\in \Bbb{N}$ to be the cardinality of
the holonomy group of $L_p$. If $d$ is the dimension of the leaves then we
define volume function of $\Cal {F}$:

$$
T:\Omega \longrightarrow \Bbb{R}^{+},
\text{ }
T(p)=n(p){\int_{L_p} }\omega ^d 
$$

\proclaim{Lemma 1}
$T${\ is a continuous locally constant function in }$\Omega $\bf{.%
}
\endproclaim

\rm
\demo{proof}
The continuity is obvious. We have to prove that $T$ is locally
constant. To do this we have just to observe that it is constant in the
residual subset of $\Omega $, formed by the union of leaves without
holonomy[G,p. 96]. By the Reeb local stability theorem there is a saturated
neighborhood for each leaf in this set where all leaves are homologous. Then
using the closedness of $\omega ^d$ and Stokes Theorem we prove the lemma.%
\enddemo

\it
Remark \rm- In fact, the proof of this lemma is essentially already contained
in [Ho].

\parindent=0pt

\vskip0.6cm

\largebf 3. A Lemma about Diff($\Bbb{C}^n,0$)
\vskip0.4cm

\rm

In 1905 Burnside[Bu] proved that if $G$ is a subgroup of $GL(n,F)$, where $F$
is a field of characteristic zero, with exponent $e$, then $G$ is finite
with $cardinality(G)\leq e^{n^3}$. Recalling that a group has exponent $e$
if every element $g$ belonging to the group is such that $g^e=1$. From the
generalization of this result by Herzog-Praeger[HP] we obtain :

\proclaim{Lemma 2}
{If }$G${\ is a subgroup of }$Diff(\Bbb{C}^n,0)${\ with
exponent }$e${\ then }$G${\ is finite with }$%
cardinality(G)\leq e^n.$
\endproclaim
\rm
\demo{proof}
If for each element of $G$ we consider its derivative we obtain a subgroup
of $GL(n,\Bbb{C)}$ with exponent $e$. Thus we only have to prove that the
normal subgroup $G_0$ of $G$, formed by its elements tangent to the identity
is the trivial group.

Let $g\in G_0$, then $g^e=Id$. Defining $H(x)={\sum_{i=1}^{e}}Dg(0)^{-i}g^i(x)$, we see that :
$$
H\circ g(x)=Dg(0)Dg(0)^{-1}
{\sum_{i=1}^{e} Dg(0)^{-i}g^{i+1}(x)=Dg(0)H(x)}
$$
Hence g is conjugated to its linear part, and therefore g must be the
identity.
\enddemo

\pagebreak

\vskip0.6cm

\largebf 4. Proof of the Theorem 1

\vskip0.4cm

\rm
Let $\Cal {F}$ be as in the theorem. Consider the connected component of $%
\Omega $ containing the leaf $L$ that is compact and with finite holonomy,
and call it $\Omega _L$. The volume function $T$ is constant in $\Omega _L$
by Lemma 1, so if $p\in \partial \Omega _L$ we have that the leaf through $p$
is aproximated by leaves with uniformly bounded volume, so it has bounded
volume and is compact(here we use the fact that the manifold is compact to
achieve the compactness of the leaf). The holonomy group of $\Omega _L$ has
finite exponent, because for any transversal $\Sigma $ of $L_p$, $\Sigma
\cap \Omega _L$ will be an open set such that every leaf of $\Omega _L$ cuts
it in at most $m$ points. Thus for every holomy germ $h$ of $L_p$,\thinspace $%
(h^{m!})_{|\Sigma \cap \Omega _L}=Id$. Analytic continuation implies that $%
h^{m!}=Id$. Using Lemma 2, we see that $\partial \Omega _L=\emptyset $, and
prove the theorem.

The Corollary follows observing that the set of leaves without holonomy is
residual and that we don't need the compactness of the manifold to assure
that a limit leaf is compact.\ Then the holonomy group of each leaf is finite
and by the results of Epstein[Ep] we get the consequences.%

\vskip0.2cm

\it Remark - \rm The same proof works in a more general
context. We have just to suppose that our foliation is transversely
quasi-analytic and that there is a closed form which is positive on the
(n-q)-planes of the distribution associated to the foliation.

\vskip0.5cm
 
\bf References

\vskip0.4cm

\smallrm

[Br] Brunella, M. :A global stability theorem for transversely holomorphic
foliations, Annals of Global Analysis and Geometry 15(1997), 179-186

[Bu] Burnside, W. :\ Proc. London Math. Soc. (2) 3 (1905), 435-440

[Ca-LN] Camacho, C. and Lins Neto, A. : Geometric theory of foliations, Birkhauser, 1985

[EMS] Edwards, R., Millett K. and Sullivan D. : Foliations with all leaves
compact, Topology 16(1977), 13-32

[Ep] Epstein, D.B.A. :\ Foliations with all leaves compact, Ann. Inst.
Fourier 26, 1(1976), 265-282

[G] Godbillon, C. : Feuilletages, \'{e}tudes g\'{e}om\'{e}triques,
Birkh\"{a}user, Basel, 1991

[HP] Herzog M. and Praeger C. : On the order of linear groups with fixed
finite exponent, Jr. of Algebra 43(1976), 216-220

[Ho] Holmann, H. :\ On the stability of holomorphic foliations, LNM
798(1980), 192-202

[Su] Sullivan, Dennis :\ A counterexample to the periodic orbit conjecture,
Inst. Hautes \'{E}tudes Sci. Publ. Math. 46(1976),5-14

\vskip0.8cm

\smallrm 
Jorge Vit\'{o}rio Bacellar dos Santos Pereira 

Email : jvp\@impa.br 

Instituto de Matemática Pura e Aplicada, IMPA 

Estrada Dona Castorina, 110 - Jardim Botânico

22460-320 - Rio de Janeiro, RJ, Brasil

\enddocument